\documentclass{amsart}

\usepackage{amscd}
\usepackage{amsmath}
\usepackage{amssymb}
\usepackage{amsthm}
\usepackage{epsf}
\usepackage{latexsym}
\usepackage{verbatim}
\input epsf.tex
\newtheorem{theorem}{Theorem}[section]
\newtheorem{definition}[theorem]{Definition}
\newtheorem{lemma}[theorem]{Lemma}
\newtheorem{corollary}[theorem]{Corollary}

\newtheorem{proposition}[theorem]{Proposition}
\theoremstyle{definition}
\newtheorem{example}[theorem]{Example}

\begin{document} 

\bibliographystyle{plain} 

\title{Measurable Sensitivity}

\author[J. James]{Jennifer James}
\address[Jennifer James]{ Williams College, MA 01267, USA}
\email{Jennifer.E.James@williams.edu}
\author[T. Koberda]{Thomas Koberda}
\address[Thomas Koberda]{  University of Chicago, Chicago, IL 60637, USA }
\email{ koberdat@uchicago.edu}
\author[K. Lindsey]{Kathryn Lindsey}
\address[Kathryn Lindsey]{ Williams College, MA 01267, USA}
\email{Kathryn.A.Lindsey@williams.edu}
\author[C.E. Silva]{Cesar E. Silva}
\address[Cesar Silva]{Department of Mathematics\\
     Williams College \\ Williamstown, MA 01267, USA}
\email{csilva@williams.edu}
\author[P. Speh]{Peter Speh}
\address[Peter Speh]{ Princeton University, Princeton, NJ 08544, USA}
\email{pspeh@princeton.edu} 
\date{\today}

\begin{abstract} We introduce the notion of measurable sensitivity, a
measure-theoretic version of the condition of sensitive dependence
on initial conditions.  It is a consequence of light mixing, implies
a transformation has only finitely many eigenvalues, and does not
exist in the infinite measure-preserving case.  Unlike the
traditional notion of sensitive dependence, measurable sensitivity
carries up to measure-theoretic isomorphism, thus ignoring the
behavior of the function on null sets and eliminating dependence on
the choice of metric.
\end{abstract}

\maketitle 


\section{Introduction}

Sensitive dependence on initial conditions, which is widely
understood to be one of the central ideas of chaos, is a
topological, rather than measurable notion.  It was introduced by
Guckenheimer in \cite{Guckenheimer}. A transformation $T$ on a
metric space $(X,d)$ is said to exhibit \textbf{sensitive
dependence} with respect to $d$ if there exists a $\delta> 0$ such
that for all $\varepsilon > 0$ and for all $x \in X$, there exists an
$n \in \mathbb{N}$ and a $y \in B_{\varepsilon}(x)$ such that
$d(T^n(x),T^n(y)) > \delta$. The notion of sensitive dependence has been
studied extensively and the reader is referred to
\cite{BBCDS}, \cite{GlasnerWeiss} and \cite{AEAJBK} for recent results. The  relationship 
between measure theoretic notions, such as weak mixing, and
sensitive dependence is studied in  \cite{GlasnerWeiss}, \cite{AEAJBK}, \cite{BGK},  and 
\cite{HYW}.  A stronger
notion of sensitivity, called strong sensitivity, was introduced by
\cite{AbrahamBC}.  A transformation $T$ on a metric space $(X,d)$ is
said to exhibit \textbf{strong sensitive dependence} with respect to
$d$ if there exists a $\delta> 0$ such that for all $\varepsilon
> 0$ and for all $x \in X$, there exists an $N \in \mathbb{N}$ so
that for all integers $n \geq N$ there exists a $y \in
B_{\varepsilon}(x)$ such that $d(T^n(x),T^n(y))
> \delta$. Both sensitive dependence and strong sensitivity are
topological notions, depending on both the choice of metric and the
behavior of the transformation on null sets. We introduce the
notions of measurable sensitivity and weak measurable sensitivity,
which are ergodic-theoretic versions of strong sensitive dependence
and sensitive dependence, respectively.  

In Section~\ref{S:2} we show that
a doubly ergodic (a condition equivalent to weak mixing for finite measure-preserving
transformations) nonsingular transformation is weak measurably sensitive, 
that a lightly mixing 
nonsingular transformation (for example, a mixing finite measure-preserving transformation) is measurably sensitive, and that measurable sensitivity does
not imply weak mixing.  In Section~\ref{S:3} we show that if an ergodic nonsingular transformation
$T$ is measurably sensitive, then there exists an integer $n>0$ so that $T^n$ has $n$ invariant 
subsets, and the restriction of $T^n$ to each of these subsets is weakly mixing.  Section~\ref{S:4} shows 
that if an ergodic finite measure-preserving transformation $T$ is measurably sensitive, then
there is an integer $n$ so that $T^n$ has $n$ invariant sets of positive measure covering $X$ a.e. and
such that the restriction of $T^n$ to each is lightly mixing.  The final section shows that an 
ergodic infinite measure-preserving transformation cannot be measurably sensitive (though it can be
weak measurably sensitive by Section~\ref{S:2}).

All of our spaces are Lebesgue spaces with a probability or a $\sigma$-finite measure defined on them.
We assume the measures to be regular.

Throughout the paper, $(X,S(X),\mu)$ denotes a  Lebesgue space $X$
with a positive, finite or $\sigma$-finite non-atomic measure $\mu$,
and $S(X)$ the collection of $\mu$-measurable subsets of $X$. It is
a standard fact that any two such spaces are isomorphic under a
nonsingular isomorphism.   We let $d$ denote a metric on $X$. We
shall say that a metric is \textbf{good} if all nonempty open sets
have positive measure. When $X$ has a good metric we assume the
measures defined on $X$ are regular.  For convenience, given two
non-empty sets $A$ and $B$ in a space $X$ with metric $d$, define
$d(A,B) = \inf \{d(a,b):a \in A, b \in B\}$.

\begin{definition} \label{t:msdef}
Let $(X,\mu,T)$ be a nonsingular  dynamical system.  $T$ is
\textbf{measurably sensitive} if whenever a dynamical system
$(X_1,\mu_1,T_1)$ is measure-theoretically isomorphic to $(X,\mu,T)$
and $d$ is a good metric on $X_1$, then there exists  a $\delta
> 0$ such that for all $x \in X_1$ and all $\varepsilon > 0$ there
exists  an $N \in \mathbb{N}$ such that for all integers $n \geq N$
\[\mu_1 \{y\in B_{\varepsilon}(x): d(T_1^nx,T_1^ny) > \delta \} > 0.\]
\end{definition}

\begin{definition} \label{t:wmsdef}
Let $(X,\mu,T)$ be a nonsingular dynamical system.  $T$ exhibits
\textbf{weak measurable sensitivity} if whenever a dynamical system
$(X_1,\mu_1,T_1)$ is measure-theoretically isomorphic to $(X,\mu,T)$
and $d$ is a good metric on $X_1$ then there exists $\delta
> 0$ such that for all $x \in X$ and $\varepsilon > 0$ there exists an
$n \in \mathbb{N}$ such that $\mu_1 \{y\in B_{\varepsilon}(x):
d(T_1^n(x),T_1^n(y)) > \delta \} > 0$.
\end{definition}

The real number $\delta$ in the above definition is referred to as a
sensitivity constant.

\begin{example} \label{t:example}
Consider a measure space X consisting of two copies of the circle
$S^1$, labeled $S_1$ and $S_2$.  Define a metric $d$ on $X$ as
follows: for two points $x \in S_i$ and $y\in S_j$, if $i=j$ define
$d(x,y)$ to be the the minimal arclength between points $x$ and $y$, and
if $i \neq j$ we let $d(x,y) = r$ for some fixed $r> \frac{\pi}{2}$.
In each copy of $S^1$, pick out the orbit of a fixed point $z$ under
a fixed irrational rotation $R$ on $S^1$, and denote this set by
$M$. Define $T:X \rightarrow X$ as follows: for $x\in M$, let $T(x)
= x$, and for $x\in X\setminus M$, let $T(x)$ map $x$ to $R(x)$ in
the \textit{other} copy of $S^1$.

Because the points in $M$ go ``far away" from those in $M^c$, it is
clear that the system exhibits sensitive dependence, but without
this null set sensitive dependence would fail.  Thus, the system
does not exhibit measurable or weak measurable sensitivity. In fact,
$T^n(B_{\varepsilon}(x)) = B_{\varepsilon}(R^2(x))$ for all $n = 0
\pmod{2}$.
\end{example}

\begin{proposition}
Let $X$ be an interval of finite length in $\mathbb{R}$ and let $d$
be the standard Euclidean metric on $X$. If a continuous
transformation $T:X\rightarrow X$ is sensitive with respect to $d$,
then $T$ is strongly sensitive with respect to $d$.
\end{proposition}

\begin{proof} Suppose $T$ is sensitive with sensitivity constant $\delta$.
Let $I_1$,...,$I_n$ be disjoint (except at endpoints) intervals with
closed or open endpoints which cover $X$ and each have length
shorter than $\frac{\delta}{2}$. Every interval of length at least
$\delta$ must contain one of these intervals.  Since $T$ is
sensitive, for each $1 \le j \le n$ there must exist a natural
number $m_j$ so that $T^{m_j}(I)$ has length at least $\delta$.
Consequently, for any interval $I$ with length at least $\delta$,
and any $n\in \mathbb{N}$, $T^{n}(I)$ contains one of $T^{k}(I_i)$
where $1\le i \le n$ and $1\le k \le m_i$.  Let $\delta'$ be one
third of the minimum of the lengths of these intervals.  Let $x\in
X$ and $\varepsilon > 0$. Since $T$ is sensitive, there must be some
$k_0$ so that $T^{k_0}(B_\varepsilon(x))$ contains an interval of
length $\delta$ and hence for any $k > k_0$,
$T^{k}(B_\varepsilon(x))$ contains an interval of the length at
least $3\delta'$ and hence contains a point whose distance is at
least $\delta'$ from $T^{k}(x)$.  Hence, $T$ is strongly sensitive
with strong sensitivity constant $\delta'$. \end{proof}

\subsection{Acknowledgements}
This paper is based on research by the Ergodic Theory group of the 2006 SMALL  summer research project at Williams College.  Support for the project was provided by National Science Foundation REU Grant DMS - 0353634 and the Bronfman Science Center of Williams College.

\section{Mixing Notions and Measurable and Weak Measurable Sensitivity}\label{S:2}
We begin by proving that double ergodicity implies weak measurable
sensitivity.  A nonsingular transformation is said to be \textbf{doubly ergodic}
if for all sets $A$ and $B$ of positive measure there exists an
integer $n > 0$ such that $\mu(T^{-n}(A) \cap A) > 0$ and
$\mu(T^{-n}(A) \cap B) > 0$.  Double ergodicity is equivalent to
weak mixing for measure-preserving transformations on finite measure
spaces~\cite{Furstenberg}, but strictly stronger than weak mixing in the
infinite measure-preserving case \cite{BFMS}.

\begin{proposition} \label{t:doubleergodicity}
If $(X,\mu,T)$ is a nonsingular, doubly ergodic dynamical system,
then $T$ exhibits weak measurable sensitivity.  In particular,
 weakly mixing, finite measure-preserving transformations
 exhibit weak measurable sensitivity.
\end{proposition}

\begin{proof}
Let $(X_1,\mu_1,T_1)$ be measure-theoretically isomorphic to
$(X,\mu,T)$, and let $d$ be a good metric on $X_1$. By the
definition of a good metric, there exist sets $A$ and $C$ of
positive measure in $S(X_1)$ such that $d(A,C) > 0$. Pick $\delta <
d(A,C)/2$. Let $x \in A$ and $\varepsilon
> 0$.  As $(X_1,\mu_1,T_1)$ is  doubly ergodic,  there exists an $n
\in \mathbb{N}$ such that $\mu (T_1^{-n}C \cap B_{\varepsilon}(x))
> 0$ and $\mu (T_1^{-n}A \cap B_{\varepsilon}(x))
> 0$.  Therefore, $\mu_1 \{y \in B_{\varepsilon}(x): T_1^n(y) \in A\} >
0$ and $\mu \{y \in B_{\varepsilon}(x):T_1^n(y) \in C\} > 0$. Since
$T_1^n(x)$ cannot be within $\delta$ of both $A$ and $C$, $\mu_1 \{y
\in B_{\varepsilon}(x) : d(T_1^n(x),T_1^n(y)) > \delta \} > 0$.
\end{proof}

We now turn our attention to measurable sensitivity and light
mixing. Recall that a system $(X,\mu,T)$ on a finite measure space
is said to be \textbf{lightly mixing} if $\liminf_{n \rightarrow
\infty} \mu(T^{-n}(A) \cap B) > 0$ for all sets $A$ and $B$ of
positive measure.

\begin{proposition}
If $(X,\mu,T)$ is a nonsingular, lightly mixing dynamical system,
then $T$ is measurably sensitive.
\end{proposition}

\begin{proof}
Let $(X_1,\mu_1,T_1)$ be measure-theoretically isomorphic to
$(X,\mu,T)$, and let $d$ be a good metric on $X_1$. By the
definition of a good metric, there exist sets $A, C \subseteq X_1$
of positive measure with $d(A,C) > 0$. Pick $\delta <
\frac{d(A,C)}{2}$. Let $x \in X_1$ and $\varepsilon> 0$. By light
mixing, there exists an $N \in \mathbb{N}$ such that for any integer
$n \geq N$, $\mu_1(T_1^{-n}(C) \cap B_{\varepsilon}(x))> 0 \
\textrm{and} \ \mu_1(T_1^{-n}(A) \cap B_{\varepsilon}(x))> 0.$ Thus,
$n \geq N$ implies $\mu_1\{y \in B_{\varepsilon}(x):T_1^n(y) \in A\}
> 0$ and $\mu_1\{y \in B_{\varepsilon}(x):T_1^n(y) \in C\} > 0$.
Since $d(A,C)> 2\delta$, $T_1^n(x)$ cannot be within $\delta$ of
both $A$ and $C$.  As $T_1^n(B_{\varepsilon}(x))$ intersects both $A$
and $C$ in sets of positive measure, $\mu_1\{y \in
B_{\varepsilon}(x) : d(T_1^n(x),T_1^n(y))
> \delta\}> 0$, so $T$ is measurably sensitive.
\end{proof}

In Proposition \ref{T: weak implies light}, we prove that a finite
measure-preserving, weakly mixing transformation that is not lightly
mixing is not measurably sensitive (but is weak measurably sensitive).

\begin{lemma}
If $T_1$ is a finite measure-preserving lightly mixing
transformation on $(X,\mu)$ and $T_2$ is a rotation on two points,
then $T \times S$ is measurably sensitive and ergodic but not weakly
mixing.
\end{lemma}

\begin{proof}
It is well known that all powers of $T$ are lightly mixing.  We
claim that $T\times S$ is measurably sensitive.  $(T_1 \times
T_2)^2$ acts as a lightly mixing transformation on $X\times\{1\}$
and $X\times\{2\}$.  Let $(S,Y)$ be isomorphic to $(T_1\times
T_2,X\times\{1,2\})$ and let $g:X\times\{1,2\}\rightarrow Y$ be the
corresponding isomorphism of measure spaces. Under a good metric on
$Y$, there exist sets $A_i$ and $B_i$ in $g(X\times \{i\})$ with a
positive distance between $A_i$ and $B_i$. Any ball $B_\epsilon(y)$
will intersect at least one of $g(X\times \{1\})$ and $g(X\times
\{2\})$ with positive measure. Then, as a consequence of light
mixing, for $n$ sufficiently large, either both $S^{-n}(A_1)$ and
$S^{-n}(B_1)$ or both $S^{-n}(A_2)$ and $S^{-n}(B_2)$ intersect
$B_\epsilon(x)$ with positive measure. Consequently, $S$ exhibits
strong sensitive dependence with a sensitivity constant $\delta =
\min\{d(A_1,B_1\}, d(A_2, B_2)\}$. Hence, $S$ is strongly sensitive.
As $T_2$ is ergodic and $T_1$ is weakly mixing, $T_1 \times T_2$ is
ergodic. Finally, $T_1\times T_2$ is not weakly mixing since $-1\in
e(T_1\times T_2)$, the eigenvalue group of $T_1\times T_2$.
\end{proof}

\section{Measurable Sensitivity and Eigenvalues}\label{S:3}

We now show that if an ergodic  nonsingular transformation is measurably
sensitive, then it can have only finitely many eigenvalues.  (Recall that 
$\lambda$ is an ($L^\infty$) eigenvalue of $T$ if there is a nonzero a.e. $f\in L^\infty$ such that
$f\circ T=\lambda f$ a.e.  Also, if $T$ is ergodic and finite measure-preserving, its $L^2$ eigenfunctions are in $L^\infty$.)  This is
used to give a further characterization of measurably sensitive
transformations.  All  ($L^\infty$) eigenvalues of    ergodic transformations lie on
the unit circle.  We refer to an eigenvalue as rational if it is of
finite order and irrational if it is not.

\begin{lemma} \label{T: intervals}
Suppose an ergodic nonsingular transformation $T: X\rightarrow X$ has an
eigenfunction $f$ with an eigenvalue that is of the form $\exp(2\pi
iq)$ with $q$ irrational, with $|f|=1$. Then for any measurable set $A \subset
S^1$ of positive Haar measure, the backwards orbit of the set $f^{-1}(A)$
equals $X \mod \mu$.
\end{lemma}
\begin{proof}
 Define $h:S^1 \rightarrow S^1$ by $h(x)
= xe^{-2\pi iq}$. The pushed measure $\mu\circ f^{-1}$ that is invariant under
$h$ must be Haar measure as $h$ is an irrational rotation.  Then

\begin{eqnarray*}
\bigcup_{n=0}^{\infty}T^{-n}(f^{-1}(A)) =
\bigcup_{n=0}^{\infty}f^{-1}(h^{n}(A)) = f^{-1}(S^1) = X.
\end{eqnarray*}
\end{proof}

\begin{lemma} \label{T: irrational eval}
Suppose an ergodic nonsingular transformation $T: X\rightarrow X$ has an
eigenfunction $f$ with an eigenvalue that is of the form $\exp(2\pi
iq)$ with $q$ irrational.  Then T is not measurably sensitive.
\end{lemma}

\begin{proof}
Assume that $|f|=1$. Construct a nonsingular isomorphism from
$(X,\mu)$ to $[0,1)$ as follows.  As a consequence of Lemma \ref{T:
intervals} and countable subadditivity, each of the sets
$f^{-1}(\exp(i(2^{1-n}\pi,2^{2-n}\pi]))$ has positive measure. Then
for $n\ge 1$, there exist nonsingular isomorphisms $g_{n}$ from
$f^{-1}(\exp(i(2^{1-n}\pi,2^{2-n}\pi]))$ to $(2^{-n},2^{-n+1})$ with
Lebesgue measure. Let $g$ be the point map associated with $g_n$ on
$f^{-1}(\exp(i(2^{1-n}\pi,2^{2-n}\pi]))$ for each positive $n$. We
now define a new metric $d$ on $X$ given by $d(x,y) = |g(x)-g(y)|$.
As each isomorphism is nonsingular and open sets have positive
measure under Lebesgue measure, $d$ is a good metric. Let $x \in X$
be a point so that $|g(x) - 2^{-n}| < 2^{-n-1}$ and let $\varepsilon <
2{^{-n-1}}$. Then $B^{d}_{\varepsilon}(x) \subset
f^{-1}(\exp(i(0,2^{1-n}\pi]))$. Let $h$ be a translation by $q$ on
$\mathbb{R}/\mathbb{Z}$. Then there is a sequence
$\{m_k\}_{k=1}^{\infty}$ for which

$$h^{m_k}((0,2^{1-n}])\subset (0,2^{2-n}].$$

Consequently, $T^{m_k}(B_{\varepsilon}(x))$ is contained in
$f^{-1}(\exp(i(0,2^{2-n}\pi]))$, and so for almost every point $y\in
T^{m_k}(B_{\varepsilon}(x))$, $d(T^{m_k}(y), T^{m_k}(x)) < 2^{2-n}$.
As $n$ tends to infinity, $2^{2-n}$ tends to zero, and so there is
no possible sensitivity constant.  Hence, $T$ is not measurably
sensitive.
\end{proof}

\begin{lemma} \label{t:rationaleigenvalues}
If an ergodic  nonsingular transformation $T$ on a Lebesgue space $X$ has
infinitely many rational eigenvalues, then $T$ is not measurably
sensitive.
\end{lemma}

\begin{proof}
Since the eigenvalues of a transformation form a multiplicative
group, there exists an increasing sequence of positive integers so
that $e^{\frac{2\pi i}{a_i}}$ is an eigenvalue and $a_i|a_{i+1}$.
These eigenvalues indicate the existence of $a_i$ invariant sets
under $T^{a_i}$. There are $\frac{a_{i+1}}{a_i}$ invariant sets for
$T^{a_i+1}$ contained in each invariant set for $T^{a_i}$. The
invariant sets can be enumerated so that the invariant sets for
$T^{a_i}$ are represented by integers $A_{i,0}, ..., A_{i,a_i-1}$
and the $T^{a_{i+1}}$-invariant sets  $A_{i+1,m}$ contained in a
given $T^{a_i}$-invariant set $A_{i,n}$ have $m = 2 \pmod{a_i}$.
Then any sequence of integers $\{b_i\}_{i=1}^{\infty}$ for which $0
\le b_i < \frac{a_{i+1}}{a_i}$ corresponds to a sequence of sets
$C_{i}$ such that $C_{i+1}\subset C_{i}$, where
$C_{i}=A_{i,\sum_{j=1}^{i}b_j a_j}$. For any such sequence, the set
$\bigcap_{i=1}^{\infty}C_i$ cannot intersect any analogous set
corresponding to another sequence. There are uncountably many
possible sequences, so there must be a sequence of sets
$\{a_i\}_{i=0}^{\infty}$ so that $\bigcap_{i=1}^{\infty}C_i$ has
measure $0$ for the corresponding sequence $a_i$. Choose such a
sequence. Then letting $C_{0} = X$, the space may be expressed as

$$X = \left(\bigsqcup_{i=0}^{\infty}C_{i}\backslash C_{i+1}\right) \mbox{ (mod $\mu$)}.$$

Let $g_{i}$ be a nonsingular isomorphism from $C_{i}\backslash
C_{i+1}$ to $(2^{-i-1},2^{-i})$ with Lebesgue measure. Let $N$ be
the backwards orbit of $\bigcap_{i=1}^{\infty}C_i$, where none of
the maps $g_{i}$ are defined.  This must have measure zero. Then let
$g: X\backslash N \rightarrow X\backslash N$ be the union of the
maps $g_i$.  Let $T'$ be the restriction of $T$ to $X\backslash N$.
Let $d$ be a metric on $X\backslash N$ defined by $d(x,y) =
|g(x)-g(y)|$.  This metric is good since each map $g_i$ is
nonsingular. Choose a point $x\in C_{i}\backslash N$ and let
$\varepsilon$ be small enough so that $B^d_{\varepsilon}(x)\subset C_{i}$.
Then for each positive integer $k$,
$(T')^{ka_i}(B^d_{\varepsilon}(x))\subset C_{i}\backslash N$.  Any two
points in $C_i\backslash N$ can have a distance of at most $2^{-i}$
between them and so any sensitivity constant $\delta$ for $T'$ would
have to be at most $2^{-i}$.  Consequently, there can be no positive
sensitivity constant, and $T$ is not measurably sensitive.
\end{proof}

\begin{corollary} \label{t:finite evals}
Any   ergodic, nonsingular, measurably sensitive transformation has finitely
many eigenvalues.
\end{corollary}

\begin{proposition} \label{t:invariant sets}
If an ergodic, nonsingular transformation $T:X\rightarrow X$ is
measurably sensitive, then for some $n\in\mathbb{N}$, $T^{n}$ has
$n$ invariant subsets and the restriction of $T^{n}$ to each of
these subsets is weakly mixing.
\end{proposition}

\begin{proof}
By Corollary \ref{t:finite evals}, $T$ must have finitely many
eigenvalues. The eigenvalues form a cyclic group of finite order
$n$. Let $f$ be an eigenfunction whose eigenvalue has order $n$.
Then the sets $A_k = f^{-1}\left(\exp\left(i\left[\frac{2\pi
k}{n},\frac{2\pi (k+1)}{n}\right)\right)\right)$ for $0 \le k \le
n-1$ are invariant under $T^{n}$. First, we show that the
restriction of $T^{n}$ to each of these sets is ergodic. Suppose
there exists a $T^n$-invariant set $C$ of positive measure. Notice
that since $T^{-n}(C) = C$, then $\bigcup_{i=0}^{n-1}T^{-i}(C)$ is
$T$-invariant, and thus equals $X \mod \mu$.  Each of $T^{-i}(C)$
for $0\le i \le n-1$ must be contained in $A_k$ for a different $k$
and so $C\cap A_k = A_k \mod \mu$ for all $k$.  Hence, the
restriction of $T^{n}$ to $A_{k}$ is ergodic for $0 \le k \le n-1$.

Next, we prove that the restriction of $T^{n}$ to any of these sets
admits no eigenvalues other than $1$. Suppose $T^{n}$ restricted to
$A_0$ admits an eigenvalue $\lambda \ne 1$. For $0\le k \le n-1$,
let $h(T^{-k}(x)) = \lambda^{\frac{-k}{n}}f(x)$. Then $h$ will be an
eigenfunction of $T$ which will have order greater than $n$. Thus
the restriction of $T^{n}$ to $A_{i}$ for each $0\le i \le n-1$ is
ergodic and admits no eigenvalues other than $1$.  Thus, the
restriction must be weakly mixing.
\end{proof}

\begin{corollary}
A totally ergodic, measurably sensitive transformation is weakly
mixing.
\end{corollary}

\section{Measurable Sensitivity for Finite Measure-preserving Transformations}
\label{S:4}

In this section, we consider measurable sensitivity for
measure-preserving transformations on finite measure spaces.
Considering only such spaces, Proposition \ref{t:invariant sets} is
strengthened to include a requirement of light mixing.  We assume
spaces have total measure one.

\begin{lemma}\label{T: separate sets}
Suppose $T:X\rightarrow X$ is measure-preserving and not lightly
mixing. Then there exist sets $C_{1}$ and $D_{1}$ of positive
measure and an infinite sequence $\{n_k\}_{k=1}^{\infty}$ such that
$T^{n_k}(C_1)\cap D_{1}= \emptyset$ for all $k \in \mathbb{N}$.
\end{lemma}
\begin{proof}
From the definition of lightly mixing, we may assume that there are
sets $C$ and $D$ of positive measure so that
$\liminf_{n\rightarrow\infty}\mu(C\cap T^{-n}(D)) = 0$.  Choose an
increasing sequence of distinct natural numbers
$\{n_k\}_{k=1}^{\infty}$ so that $\mu(C\cap T^{-n_k}(D)) \le
2^{-k-1}\mu(C)$.  Let $C_1 =
C\backslash\bigcup_{k=1}^{\infty}(C\cap T^{-n_k}(D))$ and $D_{1} =
D$. Then $\mu(C_1)
> \frac{1}{2}\mu(C) > 0$ and so $T^{n_k}(C_{1})\cap D_{1} = \emptyset$ for every
$k$.
\end{proof}

The proof of the following lemma is standard and is omitted.

\begin{lemma} \label{T: increasing continuous}
Let $f:(0,1)\rightarrow(0,1)$ be a continuous function with $f(x) >
x$ for every $x\in (0,1)$.  Then for every $x\in (0,1)$,
$\lim_{n\rightarrow\infty}f^{n}(x) = 1$.
\end{lemma}

\begin{lemma}\label{T: sequences of subsets}
Let $T:X\rightarrow X$ be a finite measure-preserving, weakly mixing
and not lightly mixing transformation. Then there exist sequences of
measurable sets $\{C_i\}_{i=1}^{\infty}$ and
$\{D_i\}_{i=1}^{\infty}$ satisfying the following properties.
\begin{enumerate}
\item \label{I: positive measure}
$\mu(C_i) > 0$ and $\mu(D_i) > 0$
\item \label{I: C subset}
$C_i \subset C_{i-1}$ for $i > 1$.
\item \label{I: D subset}
$D_{i-1} \subset D_i$ for $i > 1$.
\item \label{I: D growth}
$\lim_{i\rightarrow\infty}\mu(D_i) = 1$
\item \label{I: C growth}
$\lim_{i\rightarrow\infty}\mu(C_i) = 0$
\item\label{I: subsequence}
There is a sequence $\{n_k\}_{k=1}^{\infty}$ so that
$T^{n_k}(C_i)\cap D_{i}= \emptyset$.
\end{enumerate}
\end{lemma}
\begin{proof}
We let $C_1$, $D_1$, and $\{n_k\}_{k=1}^{\infty}$ be as defined in
Lemma \ref{T: separate sets}. They clearly satisfy properties 1-3
and property 6. For the inductive step, assume that $C_i$ and $D_i$
have been chosen to satisfy these properties for all $i \leq j$. As
a consequence of weak mixing, there is a zero density subset
$E_1\subset \mathbb{N}$ such that

$$\lim_{n\rightarrow\infty, n\not\in E_1}\mu(T^{-n}(C_j)\cap
C_j) = \mu(C_j)^{2}.$$

As a result, the set of values $n$ with $\mu(T^{-n}(C_j)\cap C_j)
> 0$ has density 1. Similarly, there exists a zero
density subset $E_2\subset \mathbb{N}$ such that

$$\lim_{n\rightarrow\infty, n\not\in E_2}\mu(T^{-n}(D_{j})\cap D_j) = \mu(D_j)^2.$$
Consequently, there is a natural number $m_{j}$ so that
$\mu(T^{-m_{j}}(C_j)\cap C_j)> 0$ and $\mu(T^{-m_{j}}(D_j)\cap D_j)
< \frac{1}{2}(\mu(D_j)^2 + \mu(D_j))$. Then
$\mu(T^{-m_{j}}(D_j)\backslash D_j)
> \frac{1}{2}(\mu(D_j) - \mu(D_j)^{2})$.  For positive integers $j$, we let $C_{j+1} =
C_{j}\cap T^{-m_j}(C_j)$ and let $D_{j+1} = D_{j} \cup
T^{-m_j}(D_j)$. Properties \ref{I: C subset} and \ref{I: D subset}
are clear from the definitions of $C_{j+1}$ and $D_{j+1}$.  Property
\ref{I: positive measure} for $D_{j+1}$ follows from the fact that
$D_0$ has positive measure and $D_i \subset D_{i+1}$ for all $i \le
j$. Property \ref{I: positive measure} for $C_{j+1}$ follows from
the fact that $m_j$ is chosen to make $C_{j+1}$ have positive
measure. Property \ref{I: D growth} follows from Lemma \ref{T:
increasing continuous} applied to the function $\frac{3}{2}x -
\frac{1}{2}x^2$ and the lower bound for the measure of $D_{i+1}$ in
terms of $D_{i}$. To show Property \ref{I: subsequence}, it suffices
to see that
\begin{eqnarray*}
T^{n_k}(T^{-m_j}(C_j)\cap C_j)\cap (T^{-m_j}(D_j)\cup D_j)\\
\subset T^{-m_j}(T^{n_k}(C_j)\cap D_j)\cup (T^{n_k}(C_j)\cap D_j) =
\emptyset.
\end{eqnarray*}
Property \ref{I: C growth} follows from property \ref{I: D growth}
and the fact that $T^{-n_k}(D_i)$ and $C_i$ must be disjoint.

\end{proof}

\begin{proposition} \label{T: weak implies light}
If $T:X\rightarrow X$ is weakly mixing, finite measure-preserving
and not lightly mixing, then $T$ is not measurably sensitive.
\end{proposition}

\begin{proof}
The space $X$ can be decomposed as $$X = \left(\bigsqcup_{i,j =
0}^{\infty} (C_{i}\backslash C_{i+1})\cap(D_{j+1}\backslash
D_{j})\right)\mod \mu$$ where $C_{0} = X$ and $D_{0} = \emptyset$
and $C_i$ and $D_i$ are as in Lemma~\ref{T: sequences of subsets}.
Let $g_{i,j}$ be a nonsingular isomorphism from $(C_{i}\backslash
C_{i+1})\cap(D_{j+1}\backslash D_{j})$ to $(2^{-j} +
2^{-i-1-j},2^{-j} + 2^{-i-j})$ with Lebesgue measure whenever
$(C_{i}\backslash C_{i+1})\cap(D_{j+1}\backslash D_{j})$ has
positive measure.  Let $N$ be the backwards orbit of the points
where no $g_{i,j}$ is defined. This set has measure zero, so the
restriction of $T$ to $X\backslash N$ is isomorphic to $T$.  Let
$T'$ denote this restriction. Define function $g$ on $X\backslash N$
by letting $g(x) = g_{i,j}(x)$. Then $d(x,y) = |g(x)-g(y)|$ is a
metric on $X\backslash N$.  A ball around a point $x\in X\backslash
N$ must have positive measure as each of the maps $g_{i,j}$ is
nonsingular, so the metric is good.  Note that
$g(D_{j})\subset(2^{-j},1)$ and $g(D_{j}^{c})\subset(0, 2^{-j})$.
Let $x$ be a point in $(C_{i}\backslash N)\cap D_j$ for some $j$
where such a point exists, and choose $\varepsilon$ so that
$B^{d}_{\varepsilon}(x)\subset C_{i}\cup D_j$.  From property \ref{I:
subsequence} of Lemma \ref{T: sequences of subsets}, there is a
sequence $\{n_{k}\}_{k=1}^{\infty}$ so that
$g((T')^{n_{k}}(B_{\varepsilon}(x)))\subset g^{-1}(0,2^{-i})$.  Since
any sensitivity constant $\delta$ must be smaller than $2^{-i}$ for
each positive integer $i$, there is no possible sensitivity
constant. Hence, $T'$ does not exhibit strong sensitive dependence
for any good metric $d$ and so $T$ is not measurably sensitive.
\end{proof}

\begin{lemma} \label{T: components}
Suppose a measurably sensitive transformation $T:X\rightarrow X$ has
finitely many invariant subsets of positive measure
$A_1$,$A_2$,..,$A_n$. Then the restriction of $T$ to each subset is
measurably sensitive.
\end{lemma}

\begin{proof}
Suppose the restriction of $S = T|A_i$ is not measurably sensitive.
Let $d$ be a good metric on a set $A_i^\prime$ so that
$(S^\prime,A_i^\prime)$ is measurably isomorphic to $(S,A_i)$ and
assume that $S^\prime$ does not exhibit strong sensitivity. Let $g$
be a measure-preserving isomorphism from $X\backslash A_i$ to (0,1)
with Lebesgue measure. Let $N$ be the set of points where $g$ and
the isomorphism from $A_i^\prime$ to $A_i$ are not preserved as well
as their backwards orbits. Then $N$ must be a null set. For any two
points $x,y \in X\backslash (A_i\cup N)$, extend $d$ so that $d(x,y)
= |g(x) - g(y)|$.  Let $X^\prime = X\backslash(A_i\cup N) \cup
A_i^\prime$ and let $T^\prime:X^\prime\rightarrow X^\prime$ be equal
to the restriction of $T$ on $X\backslash(A_i\cup N)$ and equal to
$S^\prime$ on $A_i^\prime$. Then $T^{\prime}$ is measurably
isomorphic to $T$. Now $d$ is extended to a metric on $X^\prime$.
For any two points $x,y \in X\backslash (A_i\cup N)$, extend $d$ so
that $d(x,y) = |g(x) - g(y)|$.  Choose a point $y_0\in A_i^\prime$.
For points $y\in A_i^\prime$ and for $x\in X\backslash (A_i\cup N)$,
let $d(y,x)=d(x,y) = 1 + d(y,y_0)$.  It is easy to verify that the
extension of $d$ is a good metric on $X^\prime$.

As $S^\prime$ does not exhibit strong sensitivity on $A_i^\prime$,
for any $\varepsilon
> 0$, there is a ball $B_r(x)$ with $x\in A_i^\prime$ so that for some
sequence $\{n_k\}$, almost every point $y\in B_r(x)\cap A_i^\prime$
satisfies $d((S^\prime)^{n_k}(y),(S^\prime)^{n_k}(x)) <
\varepsilon$. We may assume that $r<1$. Then $B_{r}(x)= B_r(x)\cap
A_i^\prime$. As a consequence, almost every point $y\in B_{r}(x)$
satisfies $d((T^\prime)^{n_k}(y),(T^\prime)^{n_k}(x)) <
\varepsilon$. Hence, $T'$ is not strongly sensitive and so $T$ is
not measurably sensitive.
\end{proof}

\begin{theorem}
Let $T$ be an ergodic transformation on a finite measure Lebesgue
space $X$. If $T$ is measurably sensitive, then there is some
positive integer $n$ so that $T^n$ has $n$ invariant sets of
positive measure which cover almost all of $X$, and the restriction
of $T^n$ to each of the sets is lightly mixing.
\end{theorem}

\begin{proof}
Suppose $T$ is measurably sensitive.  By Proposition
\ref{t:invariant sets}, there are $n$ invariant sets for $T^{n}$,
each of positive measure. Every power of a measurably sensitive
transformation is clearly measurably sensitive, so $T^n$ must be
measurably sensitive. By Lemma \ref{T: components}, the restriction
of $T^{n}$ to any of the sets must be measurably sensitive and, by
Proposition \ref{t:invariant sets}, the restriction must be weakly
mixing. Consequently, Proposition \ref{T: weak implies light}
indicates that the restriction must be lightly mixing.
\end{proof}

\section{Infinite Measure Spaces}
\label{S:5}

While the existence of finite measure-preserving, measurably
sensitive transformations is implied by the existence of lightly
mixing finite measure-preserving transformations, there is no
corresponding notion of light mixing for the infinite
measure-preserving case. In fact, we show that there are no ergodic,
infinite measure-preserving, measurably sensitive transformations.
There exist non-conservative ergodic nonsingular transformations
that are measurable sensitive; for example let $T:X\to X$ be a finite 
measure-preserving mixing transformation and define $S:X\times\mathbb{N}\to X\times\mathbb{N}$
by $S(x,n)= (Tx,n-1)$ if $n>1$ and $S(x,1)=(Tx,2)$.

\begin{proposition} \label{t: infinite}
There are no ergodic, infinite measure-preserving, measurably
sensitive transformations.
\end{proposition}

\begin{proof}
Let $T$ be an ergodic, measure-preserving transformation on a
$\sigma$-finite measure space $X$ with infinite measure.  Then $X =
\bigsqcup_{i=1}^{\infty} A_i$, where each set $A_i$ has positive
finite measure.  Let $D_i = \bigcup_{j=1}^{i}A_j$.  Then, as both
$A_{i+1}$ and $D_i$ have finite measure,
$\liminf_{n\rightarrow\infty}\mu(T^{-n}(D_i)\cap A_{i+1}) = 0$. Choose
an increasing sequence $\{n_{i,k}\}_{k=1}^{\infty}$ so that
$\mu(T^{-n_{i,k}}(D_i)\cap A_{i+1}) < 2^{-k-1}\mu(A_{i+1})$. Then let

$$C_i = A_{i+1}\backslash\left(\bigcup_{k=1}^{\infty}T^{-n_{i,k}}(D_i)\right)$$

The set $C_i$ has positive measure, is contained in $D_{i+1}$ and
$T^{n_{i,k}}(C_i)\cap D_i = \emptyset$ for every natural number $k$.
Let $g_{i}$ be a nonsingular isomorphism from $C_i$ to $(2^{-2i},
2^{-2i+1})$ with Lebesgue measure and let $h_i$ be a nonsingular
isomorphism $A_{i+1}\backslash C_i$ to $(2^{-2i+1}, 2^{-2i+2})$ with
Lebesgue measure whenever $A_{i+1}\backslash C_i$ has positive
measure.  Let $N$ be the set where none of the functions $h_i$ and
$g_i$ are defined as well as their backwards orbits.  This set must
have measure zero due to the nonsingularity of  $T$ and so $T'$, the
restriction of $T$ to $X\backslash N$, is measurably isomorphic to
$T$. Let $g: X\backslash N\rightarrow (0,1)$ be equal to whichever
of $h_i$ and $g_i$ is defined.  Note that $g(D_i)\subset
(0,2^{-2i})$. Then define metric $d$ on $X\backslash N$ by $d(x,y) =
|g(x)-g(y)|$. A ball in metric $d$ around any point in $X\backslash
N$ must have positive measure as each isomorphism is nonsingular
and so the metric is good. Let $x\in C_i\backslash N$. For
sufficiently small $\varepsilon
> 0$, $B_{\varepsilon}(x)\subset C_i$. Then
$(T')^{n_{i,k}}(B_{\varepsilon}(x))\subset D^{c}\backslash N$ and, as
any two points in $D_{i}^{c}\backslash N$ have a maximum distance of
$2^{-2i + 1}$ between them, any sensitivity constant must be at most
$2^{-2i + 1}$.  Consequently, $T'$ does not exhibit strong sensitive
dependence in this metric.  As the metric is good, $T$ is not
measurably sensitive.
\end{proof}

Although an infinite, ergodic, measure-preserving dynamical system
$(X,\mu,T)$ cannot be measurably sensitive, it can, however, exhibit
the desired property with respect to a good metric.

\begin{proposition}
Let $(X,\mu,T)$ be the Hajian-Kakutani Skyscraper and $d$ the
standard Euclidean metric $X$.  Then there exists $\delta
> 0$ such that for all $x \in X$ and all $\varepsilon
> 0$ there exists $N \in \mathbb{N}$ such that $\mu \{y\in
B_{\varepsilon}(x): d(T^n(x),T^n(y))
> \delta \} > 0$ for all integers $n \geq N$.
\end{proposition}

\begin{proof}
The Hajian-Kakutani Skyscraper is an infinite measure-preserving,
invertible, ergodic, rank-one transformation constructed from a
recursively defined sequence of columns consisting of left-open,
right-closed intervals. $C_0$ consists of $(0,1]$, and $C_{n+1}$ is
formed from $C_n$ by cutting $C_n$ into two equal pieces, placing
$2h_n$ spacers, which we denote $S_{n+1}$, above the right-hand half
of $C_n$, and stacking right-over-left.  The number of levels in
column $C_n$ is denoted $h_n$. Refer to \cite{HK} for a full
description of the construction.

We begin with some remarks about the structure of column $C_k$, for
$k \geq 1$. $C_k$  consists of 4 subcolumns of height $h_{k-1}$; the
bottom two subcolumns are points in $C_{k-1}$, and the top two are
points in $S_k$. $S_{k+1}$ may be thought of as consisting of 8
subcolumns of height $h_{k-1}$ and width half that of those in $C_k$
positioned above the right half of $C_k$.  For convenience, we refer
to these subcolumns in order from bottom to top by \textbf{$K_1$} to
\textbf{$K_{12}$}.

Let $\delta < \frac{3}{8}$, $x \in X$ and $\varepsilon
> 0$.  There exists an integer $\alpha > 1$ such that
$B_{\varepsilon}(x)$ contains a level, which we denote $I$, of
$C_{\alpha}$.  Define a sequence $\{l_i\}_{i=\alpha}^{\infty}$,
where $l_i$ is defined to be the first time for which the top level
of column $C_i$ is contained in $T^{l_i}(I)$. Set $N = l_{\alpha}$.
Then for an integer $n > N$, there exists a unique integer $k \geq
1$ such that $l_k < n \leq l_{k+1}$.

Denote the left half of $I$ by $L$ and the right half by $R$. Since
$T^{l_k + 8h_{k-1}}(I)$ contains the top level of column $C_{k+1}$,
there exists a unique $j \in \{1,...,8\}$ such that $(j-1)h_{k-1} <
n-l_k \leq jh_{k-1}$. Then $T^{n}(L) \subset K_j$ and $T^{n}(R)
\subset K_{j+4}$.  When $j = 1$, $T^{n}(L) \subset K_1 \subset
C_{k-1}$ and $T^{n}(R) \subset K_{j+4} \subset S_{k+1}$, so $S_{k}$ lies
between $T^{n}(R)$ and $T^{n}(L)$ on the real line.  As $\mu(S_k) >
1$, $d(T^{n}(L), T^{n}(R)) \ge 1$.  For $2 \le j \le 8$, $K_{j+3}$
lies between $T^{n}(R)$ and $T^{n}(L)$ on the real line and as $\mu
(K_{j+3}) \ge \frac{1}{4}$, $d(T^{n}(L),T^{n}(R)) \ge \frac{1}{4}$.
As $\delta < \frac{1}{8}$, $T^{n}(x)$ cannot be within $\delta$ of
both $T^{n}(L)$ and $T^{n}(R)$.  Hence $\mu\{y\in B_{\epsilon}(x):
d(T^{n}(x), T^{n}(y)) > \delta\} > 0$.
\end{proof}

\bibliographystyle{amsalpha}
\bibliography{ErgodicTheoryBib}

\begin{thebibliography}{10}

\bibitem{AbrahamBC}
Christophe Abraham, G{\'e}rard Biau, and Beno{\^{\i}}t Cadre.
\newblock Chaotic properties of mappings on a probability space.
\newblock {\em J. Math. Anal. Appl.}, 266(2):420--431, 2002.

\bibitem{AEAJBK}
Ethan Akin, Joseph Auslander, and Kenneth Berg.
\newblock When is a transitive map chaotic?
\newblock In {\em Convergence in ergodic theory and probability (Columbus, OH,
  1993)}, volume~5 of {\em Ohio State Univ. Math. Res. Inst. Publ.}, pages
  25--40. de Gruyter, Berlin, 1996.

\bibitem{BBCDS}
J.~Banks, J.~Brooks, G.~Cairns, G.~Davis, and P.~Stacey.
\newblock On {D}evaney's definition of chaos.
\newblock {\em Amer. Math. Monthly}, 99(4):332--334, 1992.

\bibitem{BGK}
Fran{\c{c}}ois Blanchard, Eli Glasner, Sergi{\u\i} Kolyada, and Alejandro
  Maass.
\newblock On {L}i-{Y}orke pairs.
\newblock {\em J. Reine Angew. Math.}, 547:51--68, 2002.

\bibitem{BFMS}
Amie Bowles, Lukasz Fidkowski, Amy~E. Marinello, and Cesar~E. Silva.
\newblock Double ergodicity of nonsingular transformations and infinite
  measure-preserving staircase transformations.
\newblock {\em Illinois J. Math.}, 45(3):999--1019, 2001.

\bibitem{Furstenberg}
H.~Furstenberg.
\newblock {\em Recurrence in Ergodic Theory and Combinatorial Number Theory}.
\newblock Princeton Univ. Press, Princeton, N.J., 1981.

\bibitem{GlasnerWeiss}
Eli Glasner and Benjamin Weiss.
\newblock Sensitive dependence on initial conditions.
\newblock {\em Nonlinearity}, 6(6):1067--1075, 1993.

\bibitem{Guckenheimer}
John Guckenheimer.
\newblock Sensitive dependence to initial conditions for one-dimensional maps.
\newblock {\em Comm. Math. Phys.}, 70(2):133--160, 1979.

\bibitem{HK}
A.~Hajian and S.~Kakutani.
\newblock An example of an ergodic measure preserving transformation defined on
  an infinite measure space.
\newblock {\em Lecture Notes in Math.}, 160:42--55, 1970.

\bibitem{HYW}
Lianfa He, Xinhua Yan, and Lingshu Wang.
\newblock Weak-mixing implies sensitive dependence.
\newblock {\em J. Math. Anal. Appl.}, 299(1):300--304, 2004.

\end{thebibliography}

\end{document}